\documentclass[11pt]{amsart}
\usepackage{graphicx}
\usepackage[top=26mm, bottom=26mm, left=32mm, right=32mm]{geometry}

\newtheorem{thm}{Theorem}[section]

\numberwithin{equation}{section}
\numberwithin{figure}{section}
\theoremstyle{definition}
\newtheorem{remark}{Remark}
\numberwithin{remark}{section}
\newtheorem{definition}{Definition}
\numberwithin{definition}{section}

\newcommand{\lsp}{\vspace{3mm}}

\newcommand{\vtwo}[2]{\left[\begin{array}{c} #1 \\ #2 \end{array}\right]}

\newcommand{\mtwo}[4]{\left[\begin{array}{cc} #1 & #2 \\ #3 & #4  \end{array}\right]}

\begin{document}

\begin{center}
\textbf{\large Rapid factorization of structured matrices via
randomized sampling}

\lsp

\textit{\small
P.G. Martinsson, Department of Applied Mathematics, University of Colorado at Boulder}

\lsp

\begin{minipage}{118mm}
\textbf{Abstract:} Randomized sampling has recently been
demonstrated to be an efficient technique for computing approximate
low-rank factorizations of matrices for which fast methods for
computing matrix vector products are available. This paper describes
an extension of such techniques to a wider class of matrices that
are not themselves rank-deficient, but have off-diagonal blocks that
are. Such matrices arise frequently in numerical analysis and signal
processing, and there exist several methods for rapidly performing
algebraic operations (matrix-vector multiplications, matrix
factorizations, matrix inversion, \textit{etc}) on them once
low-rank approximations to all off-diagonal blocks have been
constructed. The paper demonstrates that if such a matrix can be
applied to a vector in $O(N)$ time, where the matrix is of size
$N\times N$, and if individual entries of the matrix can be computed
rapidly, then in many cases, the task of constructing approximate
low-rank factorizations for all off-diagonal blocks can be performed
in $O(N\,k^{2})$ time, where $k$ is an upper bound for the numerical
rank of the off-diagonal blocks.
\end{minipage}

\end{center}

\section{Introduction}

There has recently been much interest in the development of fast
algorithms for structured matrices of different varieties. One class
of such matrices is the class of ``Hierarchically Semi-Separable''
(HSS) matrices, \cite{2005_HSS,2007_shiv_sheng,2008_jianlin}. These
matrices are characterized by a specific type of rank deficiencies
in their off-diagonal blocks (as described in Section \ref{sec:HSS})
and arise upon the discretization of many of the integral operators
of mathematical physics, in signal processing, in algorithms for
inverting certain finite element matrices, and in many other
applications, see
\textit{e.g.}~\cite{2007_shiv_sheng,2007_gu_FEM,hudson,2007_fem_inversion}.
The common occurrence of such matrices in scientific computing
motivates the development in
\cite{1994_starr_rokhlin,2008_jianlin,m_fastdirect,2006_shiv_penrose,MR2047424}
of fast algorithms for performing operations such as matrix-vector
multiplies, matrix factorizations, matrix inversions, \textit{etc}.

There currently is little consistency in terminology in discussing
structured matrices. The property that we here refer to as the
``HSS'' property also arises under different names in a range of
other publications, for instance
\cite{m_fastdirect,1996_mich_elongated,1994_starr_rokhlin,hudson}.
It is also closely related to the ``$\mathcal{H}^{2}$-matrices''
discussed in \cite{H2_matrix,2007_borm_H2a,2007_borm_H2b}. The
methods described in the present paper are directly applicable to
the structures described in
\cite{m_fastdirect,1996_mich_elongated,hudson}, and with minor
modifications to the structures in
\cite{H2_matrix,2007_borm_H2a,2007_borm_H2b}.

The observation that many matrices that arise in scientific
computing have off-diagonal blocks that can be approximated well by
low-rank matrices underlies many ``fast'' methods such as the Fast
Multipole Method \cite{rokhlin1987,rokhlin1997}, panel clustering
\cite{hackbusch_1987}, Barnes-Hut \cite{barnes_hut},
$\mathcal{H}$-matrices \cite{hackbusch}, \textit{etc}. It is
important to note that the HSS property imposes stronger conditions
on the off-diagonal blocks than any of the algorithms listed. We
describe these conditions in detail in Section \ref{sec:HSS}, but
loosely speaking, the HSS property requires both that large blocks
\textit{directly} adjacent to the diagonal can be approximated by
low rank matrices, and that the basis functions used be ``nested'',
\textit{i.e.}~that the basis functions on one level be expressed as
linear combinations of the basis functions on the next finer level.
The benefits obtained by imposing these stronger conditions in part
derive from faster algorithms for matrix-vector multiplies
\cite{hudson}, but more importantly from the fact that they allow
other linear algebraic operations such as matrix factorizations and
inversions to be performed in $O(N)$ time
\cite{2008_jianlin,m_fastdirect,2004_borm}. However, while there
exist very fast algorithms for manipulating such matrices once the
factors in the HSS representation are given, it is less well
understood how to rapidly compute these factors in the first place.
For matrices arising from the discretization of the boundary
integral equations of mathematical physics, \cite{mdirect} describes
a technique that is $O(N)$ in two dimensions and $O(N^{3/2})$ in
three. In other environments, it is possible to use known regularity
properties of the off-diagonal blocks in conjunction with standard
interpolation techniques to obtain rough initial factorizations, and
then recompress these to obtain factorizations with close to optimal
ranks \cite{2007_borm_H2b}. An approach of this kind with an
$O(N\log N)$ or $O(N\log^{2}N)$ complexity, depending on
circumstances, is used in \cite{hudson}.

The purpose of the present paper is to describe a fast and simple
randomized technique for computing all factors in the HSS
representation of a large class of matrices. It works in any
environment in which a fast matrix-vector multiplier is available
(for instance, an implementation of the Fast Multipole Method, or
some other legacy code) and it is possible and affordable to compute
a small number of actual matrix elements.
In order to describe the cost of the algorithm precisely, we must
introduce some notation: We let $A$ be an $N\times N$ matrix whose
off-diagonal blocks have maximal rank $k$ (in the ``HSS''-sense, see
Section \ref{sec:HSS}), we let $T_{\rm mult}$ denote the time
required to perform a matrix-vector multiplication $x \mapsto A\,x$,
we let $T_{\rm rand}$ denote the cost of constructing a pseudo
random number from a normalized Gaussian distribution, we let
$T_{\rm entry}$ denote the computational cost of evaluating an
individual entry of $A$, and $T_{\rm flop}$ denote the cost of a
floating point operation. The computational cost $T_{\rm total}$ of
the algorithm then satisfies
\begin{equation}
\label{eq:Ttotal}
T_{\rm total} \sim T_{\rm mult}\times 2\,(k+10) +
T_{\rm rand} \times N\,(k+10) +
T_{\rm entry}\times 2\,N\,k +
T_{\rm flop}\times c\,N\,k^{2},
%
\end{equation}
where $c$ is a small constant. In particular, if $T_{\rm mult}$ is
$O(N)$, then the method presented here is $O(N)$ as well.

The technique described in this paper utilizes recently published
methods for computing approximate low-rank factorizations of
matrices that are based on randomized sampling
\cite{random1,2007_PNAS}. As a consequence, there is a finite
probability that the method described here may fail in any given
realization of the algorithm. This failure probability is a user
specified parameter that in principle could be balanced against
computational cost. In practice, the probability of failure can at a
low cost be made entirely negligible. To be precise, in equation
(\ref{eq:Ttotal}), the number ``10'' in the first and the second
terms on the right hand side is chosen to yield a failure
probability that is provably less than $10^{-5}$, and appears to
actually be much smaller still. We note that when the method does
not fail, the accuracy of the randomized scheme is very high; in the
environment described in this paper, relative errors of less than
$10^{-10}$ are easily obtained.


\section{Preliminaries}
\label{sec:prel}

In this section, we introduce some notation, and list a number of
known results regarding low rank factorizations, and hierarchical
factorizations of matrices.

\subsection{Notation}
Throughout the paper, we measure vectors in $\mathbb{R}^{n}$ using their
Euclidean norm, and matrices using the corresponding operator norm.

For an $m\times n$ matrix $A$, and an integer $k = 1,\,2,\,\dots,\,\min(m,n)$,
we let $\sigma_{k}(A)$ (or simply $\sigma_{k}$ when it is obvious which matrix
is being referred to) denote the $k$'th singular value of $A$. We assume that
these are ordered so that $\sigma_{1}(A) \geq \sigma_{2}(A) \geq \cdots \geq
\sigma_{\min(m,n)}(A) \geq 0$.
We say that a matrix $A$ has ``$\varepsilon$-rank'' $k$ if $\sigma_{k+1}(A) < \varepsilon$.

We use the notation of Golub and Van Loan \cite{golub} to specify
submatrices. In other words, if $A$ is an $m\times n$ matrix with
entries $a_{ij}$, and $I = [i_{1},\,i_{2},\,\dots,\,i_{k}]$ and $J =
[j_{1},\,j_{2},\,\dots,\,j_{l}]$ are two index vectors, then we let
$A(I,J)$ denote the $k\times l$ matrix
$$
A(I,J) = \left[\begin{array}{cccc}
a_{i_{1}j_{1}} & a_{i_{1}j_{2}} & \cdots & a_{i_{1}j_{l}} \\
a_{i_{2}j_{1}} & a_{i_{2}j_{2}} & \cdots & a_{i_{2}j_{l}} \\
\vdots         & \vdots         &        & \vdots         \\
a_{i_{k}j_{1}} & a_{i_{k}j_{2}} & \cdots & a_{i_{k}j_{l}}
\end{array}\right].
$$
We let the shorthand $A(I,:)$ denote the matrix $A(I,[1,\,2,\,\dots,\,n])$,
and define $A(:,J)$ analogously.

Given a set of matrices $\{X_{j}\}_{j=1}^{l}$ 
we let $\mbox{diag}(X_{1},\,X_{2},\,\dots,\,X_{l})$ denote the block diagonal
matrix
$$
\mbox{diag}(X_{1},\,X_{2},\,\dots,\,X_{l}) = \left[\begin{array}{ccccc}
X_{1} & 0 & 0 & \cdots & 0 \\
0 & X_{2} & 0 & \cdots & 0 \\
0 & 0 & X_{3} & \cdots & 0 \\
\vdots & \vdots & \vdots & & \vdots \\
0 & 0 & 0     & \cdots & X_{l}
\end{array}\right].
$$

\subsection{Low rank factorizations}
\label{sec:ID}

We say that an $m\times n$ matrix $A$ has exact rank $k$ if there
exist an $m\times k$ matrix $E$ and a $k\times n$ matrix $F$ such
that
$$
A = E\,F.
$$
In this paper, we will utilize three standard matrix factorizations.
In describing them, we let $A$ denote an $m\times n$ matrix of rank
$k$. The first is the so called ``QR'' factorization:
$$
A = Q\,R
$$
where $Q$ is an $m\times k$ matrix whose columns are orthonormal,
and $R$ is a $k\times n$ matrix with the property that a permutation
of its columns is upper triangular. The second is the ``singular
value decomposition'' (SVD):
$$
A = U\,D\,V^{\rm t},
$$
where the $m\times k$ matrix $U$ and the $n\times k$ matrix $V$ have
orthonormal columns, and the $k\times k$ matrix $D$ is diagonal. The
third factorization is the so called ``interpolatory
decomposition'':
$$
A = A(:,J)\,X,
$$
where $J$ is a vector of indices marking $k$ of the columns of $A$,
and the $k\times n$ matrix $X$ has the $k\times k$ identity matrix
has a submatrix and has the property that all its entries are
bounded by $1$ in magnitude. In other words, the interpolatory
decomposition picks $k$ columns of $A$ as a basis for the column
space of $A$ and expresses the remaining columns in terms of the
chosen ones.

The existence for all matrices of the QR factorization and the SVD
are well-known, as are techniques for computing them accurately and
stably, see \textit{e.g.}~\cite{golub}. The interpolatory
decomposition is slightly less well known but it too always exists,
and there are stable and accurate techniques for computing it, see
\textit{e.g.}~\cite{gu1996,lowrank}. (Practical algorithms for
computing the interpolatory decomposition may produce a matrix $X$
whose elements slightly exceed $1$ in magnitude.) In the pseudo code
we use to describe the methods of this paper, we refer to such
algorithms as follows:
$$
[Q,R] = \texttt{qr}(A),\qquad [U,D,V] = \texttt{svd}(A),\qquad [X,J]
= \texttt{interpolate}(A).
$$

In the applications under consideration in this paper, matrices that
arise are typically only approximately of low rank. Moreover, their
approximate ranks are generally not known \`{a} priori. As a
consequence, the algorithms will typically invoke versions of the
factorization algorithms that take the computational accuracy
$\varepsilon$ as an input parameter. For instance,
$$
[U,D,V] = \texttt{svd}(A,\,\varepsilon)
$$
results in matrices $U$, $D$, and $V$ of sizes $m\times k$, $n\times k$, and $k\times k$,  such that
$$
||U\,D\,V^{\rm t} - A|| \leq \varepsilon.
$$
In this case, the number $k$ is the $\varepsilon$-rank of $A$, and
is of course an output of the algorithm. The corresponding functions
for computing an approximate QR factorization or an interpolatory
decomposition are denoted
$$
[Q,R] = \texttt{qr}(A,\,\varepsilon),\qquad [X,J] =
\texttt{interpolate}(A,\,\varepsilon).
$$

\begin{remark}
Standard techniques for computing partial QR and interpolatory
factorizations (Gram-Schmidt, Householder, \textit{etc}) produce
results that are guaranteed in the Frobenius norm (or some related
matrix norm that can be computed via $O(m\,n)$ methods). It is
possible to modify such techniques to measure remainders in the
$l^{2}$-operator norm, but we have found no need to utilize such
techniques since in all the applications that we have studied so
far, the relevant matrix norms are excellent predictors for each
other.
\end{remark}

%

\subsection{Construction of low-rank approximations via randomized sampling}
\label{sec:random} Let $A$ be a given $m\times n$ matrix that we
know can accurately be approximated by a matrix of rank $k$ (which
we do not know), and suppose that we seek to determine a matrix $Q$
with orthonormal columns (as few as possible) such that
$$
||A - Q\,Q^{\rm t}\,A||
$$
is small. (In other words, we seek a matrix $Q$ whose columns form an approximate ON-basis for
the column space of $A$.) When we have access to a fast technique for computing
matrix vector products $x \mapsto A\,x$, this task can efficiently be solved using
randomized sampling via the following steps:
\begin{enumerate}
\item Pick an integer $l$ that is slightly larger than $k$ (the choice $l = k+10$ will turn out to be a good one).
\item Form an $n\times l$ matrix $R$ whose entries are drawn independently from a normalized Gaussian distribution.
\item Form the product $S = A\,R$.
\item Construct a matrix $Q$ whose columns form an ON-basis for the columns of $S$.
\end{enumerate}
Note that each column of the ``sample'' matrix $S$ is a random linear combination
of the columns of $A$. We would therefore expect the algorithm described to have
a high probability of producing an accurate result provided that $l$ is sufficiently
much larger than $k$. It is perhaps less obvious that this probability depends only on
the difference between $l$ and $k$ (not on $m$ or $n$, or any other properties of $A$),
and that it approaches $1$ extremely rapidly as $l-k$ increases. The details are given in
the following theorem from \cite{random1}:

\begin{thm}
\label{thm:RAND} Let $A$ be an $m\times n$ matrix, and let $l$ and
$k$ be integers such that $l \geq k$. Let $R$ be an $n\times l$
matrix whose entries are drawn independently from a normalized
Gaussian distribution. Let $Q$ be an $m\times l$ matrix whose
columns form an ON-basis for the columns of $A\,R$. Let
$\sigma_{k+1}$ denote the minimal error in approximating $A$ by a
matrix of rank $k$:
$$
\sigma_{k+1} = \min_{\mbox{rank}(B) = k}||A - B||.
$$
Then
$$
||A - Q\,Q^{\rm t}\,A||_{2} \leq 10\ \sqrt{l\,n}\ \sigma_{k+1},
$$
with probability at least
$$
1 - \varphi(l-k),
$$
where $\varphi$ is a decreasing function satisfying, for instance,
$\varphi(8) < 10^{-5}$ and $\varphi(20) < 10^{-17}$.
\end{thm}

In this paper, we always set $l = k+10$; note however that $l-k$ is
a user defined parameter that can be set to balance the risk of
failure against computational cost.

\begin{remark}
\label{remark:k_correct} In practical applications, it is in Step
(4) of the algorithm described above sufficient to construct an
ON-basis for the columns of the sample matrix $S$ that is accurate
to precision $\varepsilon$. Under very moderate conditions on decay
of the singular values of $A$, the number of basis vectors actually
constructed will be extremely close to the optimal number,
regardless of the number $l$.
\end{remark}

\begin{remark}
The approximate rank $k$ is rarely
known in advance. In a situation where a single matrix  $A$ is to be
analyzed, it is a straight-forward matter to modify the algorithm
described here to an algorithm that adaptively determines the
numerical rank by generating a sequence of samples from the column
space of $A$ and simply stopping when no more information is added.
In the application we have in mind in this paper, however, the
randomized scheme will be used in such a way that a single random
matrix will be used to create samples of a large set of different
matrices. In this case, we choose a number $l$ of random samples
that we are confident exceeds the numerical rank of all the matrices
by at least $10$. Note that the bases constructed for the various
matrices will have the correct number of elements due to the
observation described in Remark \ref{remark:k_correct}.
\end{remark}

\begin{remark}
\label{remark:id_rand}
The randomized sampling technique is particularly effective when
used in conjunction with the interpolatory decomposition. To
illustrate, let us suppose that $A$ is an $n\times n$ matrix of rank $k$
for which we can rapidly evaluate the maps $x \mapsto A\,x$ and
$x \mapsto A^{\rm t}\,x$. Using the randomized sampling
technique, we then construct matrices $S^{\rm col} = A\,R$
and $S^{\rm row} = A^{\rm t}\,R$ whose columns span the column
and the row spaces of $A$, respectively. If we
seek to construct a factorization of $A$ without using the
interpolatory decomposition, we would then orthonormalize the
columns of $S^{\rm col}$ and $S^{\rm row}$,
$$
[Q^{\rm col},\,Y^{\rm col}] = \texttt{qr}(S^{\rm col}),
\qquad\mbox{and}\qquad
[Q^{\rm row},\,Y^{\rm row}] = \texttt{qr}(S^{\rm row}),
$$
whence
\begin{equation}
\label{eq:Afac1}
A = Q^{\rm col}\,\bigl( (Q^{\rm col})^{\rm t}\,A\,Q^{\rm row}\bigr)\,Q^{\rm row}.
\end{equation}
Note that the evaluation of (\ref{eq:Afac1}) requires $k$ matrix-vector multiplies
involving the large matrix $A$ in order to compute the
$k\times k$ matrix $(Q^{\rm col})^{\rm t}\,A\,Q^{\rm row}$.
Using the interpolatory
decomposition instead, we simply determine the $k$ rows of $S^{\rm col}$ and
$S^{\rm row}$ that span their respective row spaces,
$$
[X^{\rm col},\,J^{\rm col}] = \texttt{interpolate}((S^{\rm col})^{\rm t}),
\qquad\mbox{and}\qquad
[X^{\rm row},\,J^{\rm row}] = \texttt{interpolate}((S^{\rm row})^{\rm t}).
$$
Then we immediately obtain the factorization
\begin{equation}
\label{eq:Afac2}
A = X^{\rm col}\,A(J^{\rm col},\,J^{\rm row})\,(X^{\rm row})^{\rm t}.
\end{equation}
Note that the factorization (\ref{eq:Afac2}) is obtained by simply
extracting the $k\times k$ submatrix $A(J^{\rm col},\,J^{\rm row})$ from $A$.
\end{remark}

\begin{remark}
In practical application, the entries of the random matrix $R$ are
not ``true'' random numbers, but numbers from a ``pseudo random
number generator''. Empirical experiments indicate that the
algorithm is not at all sensitive to the quality of the random
number generator.
\end{remark}

\subsection{Hierarchically Semi-Separable matrices}
\label{sec:HSS}

In this section, we define the class of ``Hierarchically
Semi-Separable'' (HSS) matrices, and introduce notation for keeping
track of various blocks of structured matrices. A more detailed
discussion of this topic can be found in
\textit{e.g.}~\cite{2008_jianlin,2007_shiv_sheng}.

\begin{figure}
\setlength{\unitlength}{1mm}
\begin{picture}(105,41)
\put(15, 0){\includegraphics[width=90mm]{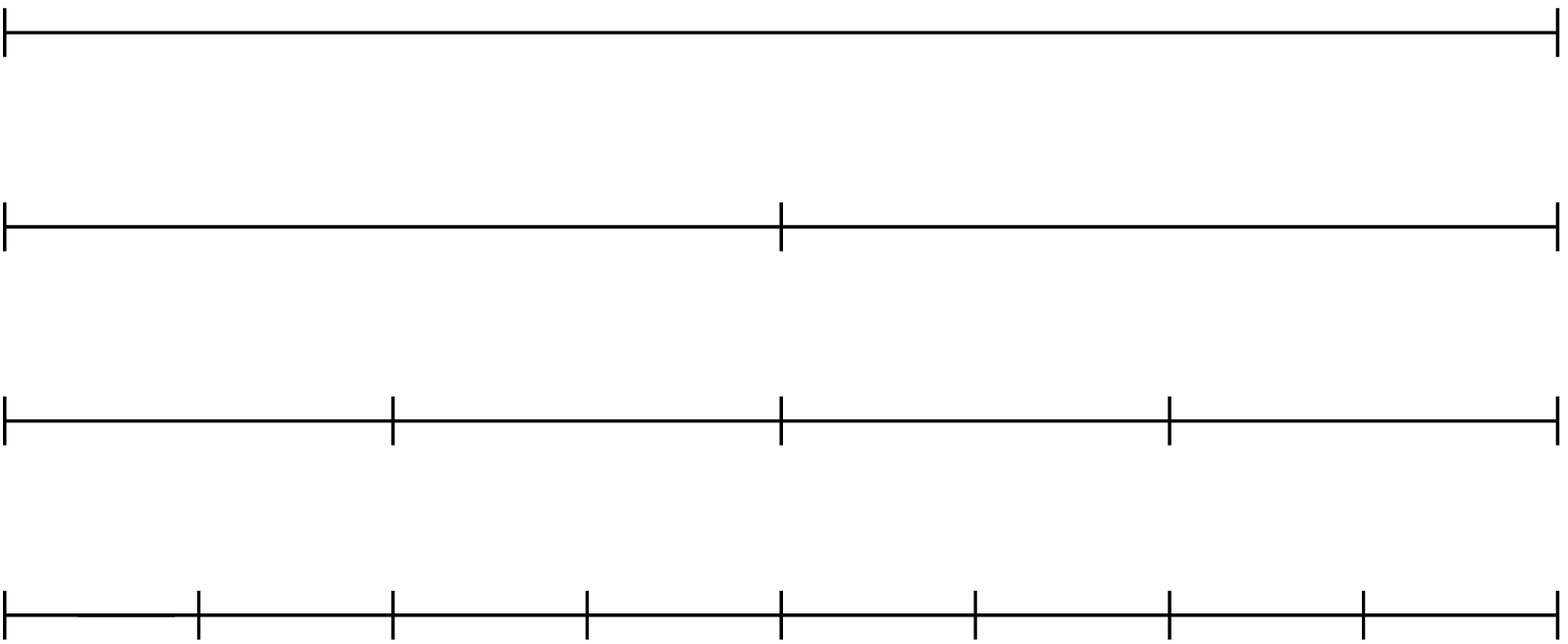}}
\put( 0, 1){Level $3$}
\put( 0,12){Level $2$}
\put( 0,23){Level $1$}
\put( 0,34){Level $0$}
\put(17, 3){$(3,1)$}
\put(28, 3){$(3,2)$}
\put(39, 3){$(3,3)$}
\put(50, 3){$(3,4)$}
\put(62, 3){$(3,5)$}
\put(73, 3){$(3,6)$}
\put(84, 3){$(3,7)$}
\put(95, 3){$(3,8)$}
\put(24,15){$(2,1)$}
\put(45,15){$(2,2)$}
\put(66,15){$(2,3)$}
\put(87,15){$(2,4)$}
\put(33,26){$(1,1)$}
\put(79,26){$(1,2)$}
\put(56,37){$(0,1)$}
\end{picture}
\caption{The binary tree of nodes.}
\label{fig:defn_tree}
\end{figure}

In order to define the HSS property for an $N\times N$ matrix $A$,
we first partition the index vector $I = [1,\,2,\,\dots\,,\,N]$ in a
hierarchy of index sets. For simplicity, we limit attention to
binary tree structures in which every level is fully populated. We
use the integers $p = 0,\,1,\,\dots,\,P$ to label the different
levels, with $P$ denoting the finest level. Thus at level $P$, we
partition $I$ into $2^{P}$ disjoint vectors
$\{I_{(P,j)}\}_{j=1}^{2^{P}}$ so that
$$
I = [I_{(P,1)},\,I_{(P,2)},\,\dots,\,I_{(P,2^{P})}].
$$
We then merge the index vectors by twos, so that for each $p =
0,\,1,\,\dots,\,P-1$, and each $j = 1,\,2,\,\dots,\,2^{p}$, we have
$$
I_{(p,j)} = [I_{(p+1,2j-1)},\,I_{(p+1,2j)}].
$$
For a given node $\tau = (p,j)$, we call the two nodes $\sigma_{1} =
(p+1,2j-1)$ and $\sigma_{2} = (p+1,2j)$ the ``children'' of $\tau$,
and we say that $\sigma_{1}$ and $\sigma_{2}$ are ``siblings''. The
tree structure is illustrated in Figure \ref{fig:defn_tree}.

\begin{figure*}
\begin{tabular}{ccccc}
\setlength{\unitlength}{1mm}
\begin{picture}(43,43)
\put(00,00){\includegraphics[height=42mm]{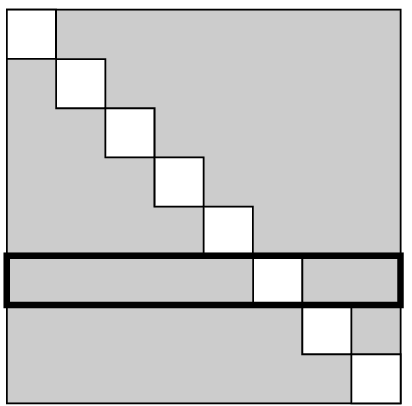}}
\end{picture}
&\hspace{5mm}&
\setlength{\unitlength}{1mm}
\begin{picture}(43,43)
\put(00,00){\includegraphics[height=42mm]{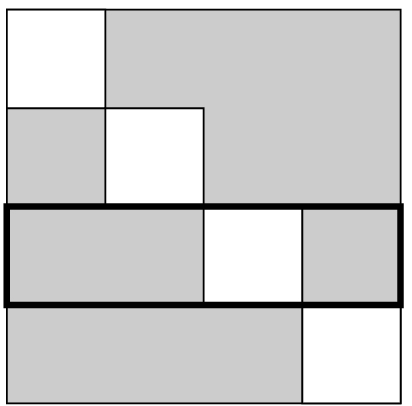}}
\end{picture}
&\hspace{5mm}&
\setlength{\unitlength}{1mm}
\begin{picture}(43,43)
\put(00,00){\includegraphics[height=42mm]{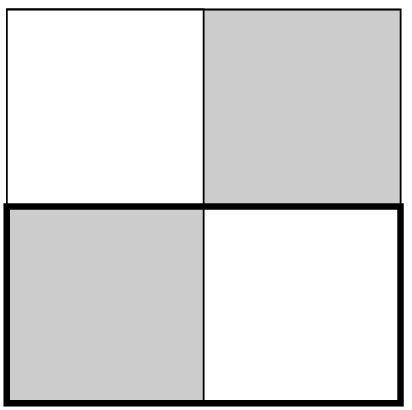}}
\end{picture}
\\
(a) && (b) && (c)
\end{tabular}
\caption{(a) The matrix $A - D^{(3)}$ with non-zero parts shaded.
The HSS row block $A_{(3,6)}^{\rm row}$ is marked with a thick border.
(b) The matrix $A - D^{(2)}$ with $A_{(2,3)}^{\rm row}$ marked.
(c) The matrix $A - D^{(1)}$ with $A_{(1,2)}^{\rm row}$ marked.}
\label{fig:HSS_def}
\end{figure*}

For any node $\tau = (p,j)$ in the tree, we define the corresponding diagonal block of $A$ via
$$
D_{\tau} = A(I_{\tau},\,I_{\tau}),
$$
and let $D^{(p)}$ denote the $N\times N$ matrix with the matrices $\{D_{(p,j)}\}_{j=1}^{2^{p}}$
as its diagonal blocks,
\begin{equation}
\label{eq:def_Dp}
D^{(p)} = \mbox{diag}(D_{(p,1)},\,D_{(p,2)},\,\dots,\,D_{(p,2^{p})}).
\end{equation}
For a node $\tau = (p,j)$, we now define the corresponding ``HSS row block''
$A_{\tau}^{\rm row}$ and ``HSS column block'' $A_{\tau}^{\rm col}$ by
$$
A^{\rm row}_{\tau} = A(I_{\tau},:) - D^{(p)}(I_{\tau},:),
\qquad\mbox{and}\qquad
A^{\rm col}_{\tau} = A(:,I_{\tau}) - D^{(p)}(:,I_{\tau}).
$$
These definitions are illustrated in Figure \ref{fig:HSS_def}.

\begin{definition}
A matrix $A$ is an ``HSS matrix'' if for some given positive integer $k$,
every HSS block has rank at most $k$.
\end{definition}

It is convenient to construct factorizations for the HSS blocks that
are ``nested'' in the sense that the bases on one level are expressed
in terms of the bases on the next finer level. To express this concept
in formulas, we suppose that for each node
$\tau$, $U_{\tau}$ denotes a matrix whose columns form
a basis for the column space of the HSS row block $A^{\rm row}_{\tau}$.
It is possible to construct the set of bases $\{U_{\tau}\}$ in such a
way that for every non-leaf node $\tau$, there exists a $2k\times k$ matrix
$\hat{U}_{\tau}$ such that
\begin{equation}
\label{eq:recursion}
U_{\tau} = \mtwo{U_{\sigma_{1}}}{0}{0}{U_{\sigma_{2}}}\,\hat{U}_{\tau},
\end{equation}
where $\sigma_{1}$ and $\sigma_{2}$ denote the two children of $\tau$.
We analogously construct a set $\{V_{\tau}\}$ of bases for the row
spaces of the HSS
column blocks $A^{\rm col}_{\tau}$ for which there exist $2k\times k$ matrices
$\hat{V}_{\tau}$ such that
$$
V_{\tau} = \mtwo{V_{\sigma_{1}}}{0}{0}{V_{\sigma_{2}}}\,\hat{V}_{\tau}.
$$

Next, let $\{\sigma_{1},\,\sigma_{2}\}$ denote a sibling pair in the tree,
and consider the offdiagonal block
$$
A_{\sigma_{1}\sigma_{2}} = A(I_{\sigma_{1}},\,I_{\sigma_{2}}).
$$
Since $A_{\sigma_{1}\sigma_{2}}$ is a submatrix of the HSS row block
$A^{\rm row}_{\sigma_{1}}$ and the HSS column block $A^{\rm col}_{\sigma_{2}}$, there
must exist a $k\times k$ matrix $B_{\sigma_{1}\sigma_{2}}$ such that
$$
A_{\sigma_{1}\sigma_{2}} =
U_{\sigma_{1}}\,B_{\sigma_{1}\sigma_{2}}\,(V_{\sigma_{2}})^{\rm t}.
$$

An HSS matrix $A$ is completely described if for every leaf node
$\tau$, we are given the matrices $D_{\tau}$ and the basis matrices
$U_{\tau}$ and $V_{\tau}$, if for all sibling pairs $\{\sigma_{1},\,\sigma_{2}\}$
we are given the matrices $B_{\sigma_{1}\sigma_{2}}$, and if for each
non-leaf node $\tau$ we are given the matrices $\hat{U}_{\tau}$ and
$\hat{V}_{\tau}$. In particular, given a vector $x$, the vector
$b = A\,x$ can be evaluated via the following steps:
\begin{enumerate}
\item For every leaf node $\tau$, calculate $\tilde{x}_{\tau} = V_{\tau}^{\rm t}\,x(I_{\tau})$.
\item Looping over all non-leaf nodes $\tau$, from finer to coarser, calculate
$\tilde{x}_{\tau} = \hat{V}_{\tau}^{\rm t}\vtwo{\tilde{x}_{\sigma_{1}}}{\tilde{x}_{\sigma_{2}}}$,
where $\sigma_{1}$ and $\sigma_{2}$ are the children of $\tau$.
\item Looping over all non-leaf nodes $\tau$, from coarser to finer, calculate\\
$\vtwo{\tilde{b}_{\sigma_{1}}}{\tilde{b}_{\sigma_{2}}} =
\mtwo{0}{B_{\sigma_{1}\sigma_{2}}}{B_{\sigma_{2}\sigma_{1}}}{0}\,
\vtwo{\tilde{x}_{\sigma_{1}}}{\tilde{x}_{\sigma_{2}}} + \hat{U}_{\tau}\,\tilde{b}_{\tau}$
(where $\tilde{b}_{\tau} = 0$ for the root node).
\item For every leaf node $\tau$, calculate $b(I_{\tau}) = U_{\tau}\,\tilde{b}_{\tau} + D_{\tau}\,x(I_{\tau})$.
\end{enumerate}
The computational cost of performing these steps is $O(N\,k)$.

\begin{remark}
The matrix $A$ can be expressed in terms of the matrices $D_{\tau}$, $U_{\tau}$,
$V_{\tau}$, and $B_{\sigma_{1}\sigma_{2}}$ as a telescoping factorization. To
demonstrate this, we introduce for each level $p = 1,\,2,\,\dots,\,P$ the
block-diagonal matrices
$$
U^{(p)} = \mbox{diag}(U_{(p,1)},\,U_{(p,2)},\,\dots,\,U_{(p,2^{p})}),
\qquad\mbox{and}\qquad
V^{(p)} = \mbox{diag}(V_{(p,1)},\,V_{(p,2)},\,\dots,\,V_{(p,2^{p})}).
$$
Moreover, we define for each non-leaf node $\tau$ the $2k\times 2k$ matrices
$$
B_{\tau} = \mtwo{0}{B_{\sigma_{1}\sigma_{2}}}{B_{\sigma_{2}\sigma_{1}}}{0}
$$
where $\sigma_{1}$ and $\sigma_{2}$ are the children of $\tau$, and set for $p = 0,\,1,\,\dots,\,P-1$
$$
B^{(p)} = \mbox{diag}(B_{(p,1)},\,B_{(p,2)},\,\dots,\,B_{(p,2^{p})}).
$$
Finally, recall from (\ref{eq:def_Dp}) that
$$
D^{(P)} = \mbox{diag}(D_{(P,1)},\,D_{(P,2)},\,\dots,\,D_{(P,2^{P})}).
$$
For simplicity, we give the factorization for the specific value $P = 3$:
\begin{equation}
\label{eq:telescope}
A = U^{(3)}\bigl(U^{(2)}\bigl(U^{(1)}\,B^{(0)}\,(V^{(1)})^{\rm t} + B^{(1)}\bigr)
(V^{(2)})^{\rm t} + B^{(2)}\bigr)(V^{(3)})^{\rm t} + D^{(3)}.
\end{equation}
The block structure of the formula (\ref{eq:telescope}) is as follows:\\
{\footnotesize
$\displaystyle\mbox{}\hspace{5mm} U^{(3)} \hspace{7mm} U^{(2)}
\hspace{3mm} U^{(1)} B^{(0)}\hspace{2mm} V^{(1)} \hspace{3mm} B^{(1)}
\hspace{8mm} V^{(2)} \hspace{8mm} B^{(2)} \hspace{20mm} V^{(3)} \hspace{26mm} D^{(3)}$
}\\
\includegraphics[width=\textwidth]{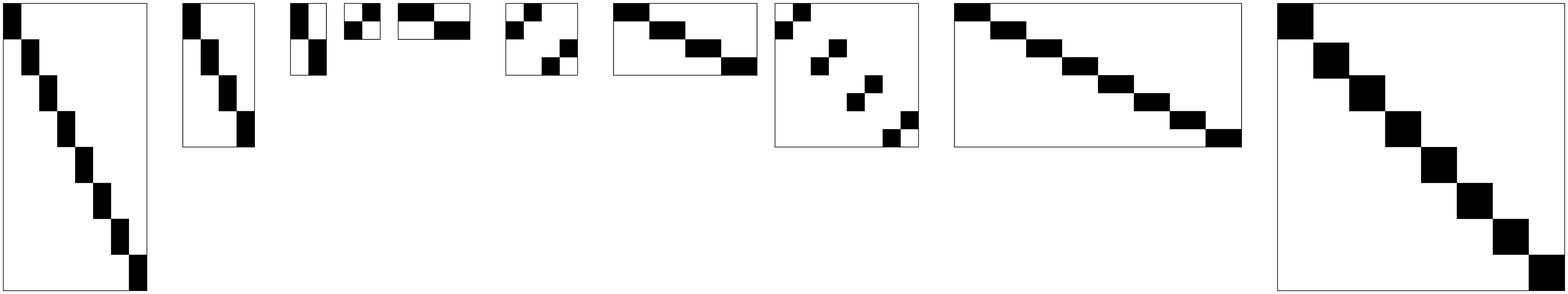}
\end{remark}

\begin{remark}
\label{remark:k}
For notational simplicity, we consider only the case where all
HSS blocks are approximated by factorizations of the same rank $k$.
In practice, it is a simple matter to implement algorithms that use
variable ranks.
\end{remark}

\begin{remark}
It is common to require the matrices $U_{\tau}$ and $V_{\tau}$
that arise in an HSS factorization of a given matrix to have
orthonormal columns. We have found it convenient to relax this
assumption and allow the use of other well-conditioned bases.
In particular, the use of interpolative decompositions
(as described in Section \ref{sec:ID}) is essential to the
performance of the fast factorization technique described in
Section \ref{sec:compHSS-ID}.
A simple algorithm for switching between the two formulations (using
orthonormal bases, or interpolatory ones) is described in Section \ref{sec:compHSS}.
\end{remark}

\section{Fast computation of HSS approximations}
\label{sec:fast}

The straight-forward way of computing the HSS factorization of matrix
would be to form all HSS blocks, and then perform dense linear algebra
operations on them to construct the required basis functions. This
approach would require at least $O(N^{2}\,k)$ algebraic operations to
factorize an $N\times N$ matrix of HSS rank $k$. In this section, we describe
how the randomized sampling techniques described in Section \ref{sec:random}
can be used to reduce this cost to $O(N\,k^{2})$.

The technique described relies crucially on the use of the interpolatory
decompositions described in Section \ref{sec:ID} in the HSS factorization.
The advantage is that the matrices $B_{\sigma_{1}\sigma_{2}}$ are then
\textit{submatrices} of the original matrix $A$ and can therefore be constructed
directly without a need for projecting the larger blocks onto the bases chosen,
\textit{cf.}~Remark \ref{remark:id_rand}.

Section \ref{sec:compHSS-ID} describes a scheme for rapidly computing
the HSS factorization of a symmetric matrix.
The scheme described in Section \ref{sec:compHSS-ID} results in a
factorization based on interpolatory bases and the
blocks $B_{\sigma_{1}\sigma_{2}}$ are submatrices of the original
matrix; Section \ref{sec:compHSS} describes how such a factorization
can be converted to one in which the bases for the HSS blocks are
orthonormal, and the blocks $B_{\sigma_{1}\sigma_{2}}$ are diagonal.
Section \ref{sec:compHSSnonsym} describes how to extend the methods
to non-symmetric matrices.

\subsection{A scheme for computing an HSS factorization of a symmetric matrix}
\label{sec:compHSS-ID}
Let $A$ be an $N\times N$ symmetric matrix that has an HSS factorization of rank $k$.
Suppose further that:
\renewcommand{\labelenumi}{(\alph{enumi})}
\begin{enumerate}
\item Matrix vector products $x \mapsto A\,x$ can be evaluated at a cost $T_{\rm mult}$.
\item Individual entries of $A$ can be evaluated at a cost $T_{\rm entry}$.
\end{enumerate}
\renewcommand{\labelenumi}{(\arabic{enumi})}
In this section, we will describe a scheme for computing an HSS factorization
of $A$ in time
\begin{equation*}
%
T_{\rm total} \sim T_{\rm mult}\times(k+10) +
T_{\rm rand} \times N\,(k+10) +
T_{\rm entry}\times 2\,N\,k +
T_{\rm flop}\times c\,N\,k^{2},
%
\end{equation*}
where $T_{\rm rand}$ is the time required to generate a pseudo random
number, $T_{\rm flop}$ is the CPU time requirement for a floating point operation and
$c$ is a small number that does not depend on $N$ or $k$.

The core idea of the method is to construct an $N\times (k+10)$ random matrix $R$,
and then construct for each level $p$, the ''sample'' matrices
$$
S^{(p)} = \bigl(A - D^{(p)}\bigr)\,R,
$$
via a procedure to be described. Then for any cell $\tau$ on level $p$,
$$
S^{(p)}(I_{\tau},:) = A^{\rm row}_{\tau}\,R,
$$
and since $A^{\rm row}_{\tau}$ has rank $k$, the columns of
$S^{(p)}(I_{\tau},:)$ span the column space of $A^{\rm row}_{\tau}$
according to Theorem \ref{thm:RAND}. We can then construct a basis
for the column space of the large matrix $A^{\rm row}_{\tau}$ by
analyzing the small matrix $S^{(p)}(I_{\tau},:)$.

What makes the procedure fast is that the sample matrices $S^{(p)}$
can be constructed by means of an $O(N)$ process from the
result of applying the \textit{entire} matrix $A$ to $R$,
$$
S = A\,R.
$$
At the finest level, $p=P$, we directly obtain $S^{(P)}$ from
$S$ by simply subtracting the contribution from the diagonal
blocks of $A$,
\begin{equation}
\label{eq:p1}
S^{(P)} = S - D^{(P)}\,R.
\end{equation}
Since $D^{(P)}$ is block diagonal with small blocks,
equation (\ref{eq:p1}) can be evaluated cheaply.
To proceed to the next coarser level, $p = P-1$, we observe that
\begin{multline}
\label{eq:p2}
S^{(P-1)} = (A - D^{(P-1)})\,R
= (A - D^{(P)})\,R - (D^{(P-1)} - D^{(P)})\,R\\
= S^{(P)} - (D^{(P-1)} - D^{(P)})\,R.
\end{multline}
Now $(D^{(P-1)} - D^{(P)})$ has only $2^{P}$ non-zero blocks.
The pattern of these blocks is illustrated for $P=3$ below:
\begin{center}
\includegraphics[width=25mm]{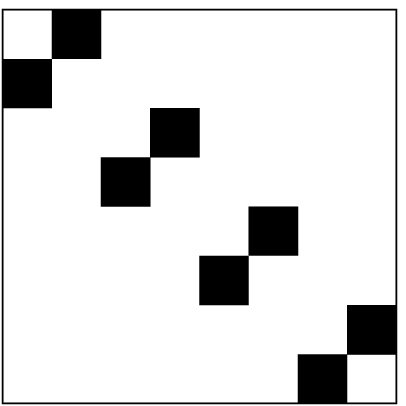}
\end{center}
Each of these blocks were compressed in the computation at level
$P$ so (\ref{eq:p2}) can also be evaluated rapidly.
The algorithm then proceeds up towards coarser levels via the formula
$$
S^{(p-1)} = S^{(p)} - (D^{(p-1)} - D^{(p)})\,R,
$$
which can be evaluated rapidly since the blocks of
$(D^{(p-1)} - D^{(p)})$ have at this point been compressed.

The condition that the bases be ``nested'' in the sense of formula
(\ref{eq:recursion}) can conveniently be enforced by using the
interpolative decompositions described in Section \ref{sec:ID}:
At the finest level, we pick $k$ rows of each HSS row block that span
its row space. At the next coarser level, we pick in each HSS row block
$k$ rows that span its row space \textit{out of the $2k$ rows that
span its two children}. By proceeding analogously throughout
the upwards pass, (\ref{eq:recursion}) will be satisfied.

The use of interpolatory decompositions has the additional benefit
that we do not need to form the entire matrices $S^{(p)}$ when
$p < P$. Instead, we work with the submatrices
formed by keeping only the rows of $S^{(p)}$ corresponding to
the spanning rows at that step.

A complete description of the methods is given with the caption ``Algorithm 1''.

\begin{figure*}
\fbox{
\begin{minipage}{\textwidth}
\begin{tabbing}
\hspace{15mm} \= \kill
\textit{Input:}  \> A fast means of computing matrix-vector products $x \mapsto A\,x$.\\
                 \> A method for computing individual entries of $A$.\\
                 \> The HSS-rank $k$ of $A$.\\
                 \> A partitioning $\{I_{(P,j)}\}_{j=1}^{2^{P}}$ of the interval $[1,\,2,\,\dots,\,N]$.\\
\textit{Output:} \> Matrices $\hat{U}_{\tau}$, $B_{\sigma_{1}\sigma_{2}}$, $D_{\tau}$ that form an
                    HSS factorization of $A$.\\
                 \> (Note that $\hat{V}_{\tau} = \hat{U}_{\tau}$ for a symmetric matrix.)
\end{tabbing}

\rule{\textwidth}{1pt}

\begin{tabbing}
\hspace{8mm} \= \hspace{8mm} \= \hspace{8mm} \= \kill
Generate an $N\times (k+10)$ Gaussian random matrix $R$.\\
Evaluate $S = A\,R$ using the fast matrix-vector multiplier.\\
\textbf{loop} over levels, finer to coarser, $p = P:(-1):1$\\
\> \textbf{loop} over all nodes $\tau$ on level $p$\\
\> \> \textbf{if} $\tau$ is a leaf node then\\
\> \> \> $I_{\rm loc} = I_{\tau}$\\
\> \> \> $R_{\rm loc} = R(I_{\tau},:)$\\
\> \> \> $S_{\rm loc} = S(I_{\tau},:) - A(I_{\tau},I_{\tau})\,R_{\rm loc}$\\
\> \> \textbf{else}\\
\> \> \> Let $\sigma_{1}$ and $\sigma_{2}$ be the two children of $\tau$.\\
\> \> \> $I_{\rm loc} = [\tilde{I}_{\sigma_{1}},\,\tilde{I}_{\sigma_{2}}]$\\
\> \> \> $R_{\rm loc} = \vtwo{R_{\sigma_{1}}}{R_{\sigma_{2}}}$\\
\> \> \> $S_{\rm loc} = \vtwo{S_{\sigma_{1}} - A(\tilde{I}_{\sigma_{1}},\tilde{I}_{\sigma_{2}})\,R_{\sigma_{2}}}
                             {S_{\sigma_{2}} - A(\tilde{I}_{\sigma_{2}},\tilde{I}_{\sigma_{1}})\,R_{\sigma_{1}}}$\\
\> \> \textbf{end if}\\
\> \> $[\hat{U}_{\tau},J_{\tau}] = \texttt{interpolate}(S_{\rm loc}^{\rm t})$\\
\> \> $R_{\tau} = \hat{U}_{\tau}^{\rm t}\,R_{\rm loc}$\\
\> \> $S_{\tau} = S_{\rm loc}(J_{\tau},:)$\\
\> \> $\tilde{I}_{\tau} = I_{\rm loc}(J_{\tau})$\\
\> \textbf{end loop}\\
\textbf{end loop}\\
For all leaf nodes $\tau$, set $D_{\tau} = A(I_{\tau},\,I_{\tau})$.\\
For all sibling pairs $\{\sigma_{1},\,\sigma_{2}\}$ set
$B_{\sigma_{1}\sigma_{2}} = A(\tilde{I}_{\sigma_{1}},\,\tilde{I}_{\sigma_{2}})$.
\end{tabbing}
\end{minipage}}

\vspace{2mm}

\textsc{Algorithm 1:} Computing the HSS factorization of a symmetric matrix.
\end{figure*}

\begin{remark}
For simplicity, Algorithm 1 is described for the case where the off-diagonal
blocks of $A$ has exact rank at most $k$, and the number $k$ is known in
advance. In actual applications, one typically is given a matrix $A$ whose
off-diagonal blocks are not necessarily rank-deficient in an exact sense,
but can to high accuracy be approximated by low-rank matrices. In this case,
Algorithm 1 needs to be modified slightly to take as an input the computational
accuracy $\varepsilon$ instead of the rank $k$, and the line
$$
[\hat{U}_{\tau},J_{\tau}] = \texttt{interpolate}(S_{\rm loc}^{\rm t})
$$
needs to be replaced by the line
$$
[\hat{U}_{\tau},\,J_{\tau}] = \texttt{interpolate}(S_{\rm loc}^{\rm t},\,\varepsilon).
$$
This directly leads to a variable rank algorithm that is typically
far more efficient than the fixed rank algorithm described.
\end{remark}

\subsection{Recompression into orthonormal basis functions}
\label{sec:compHSS}

In this section, we describe a simple post-processing step that converts
the HSS factorization resulting from Algorithm 1 (which is based on
interpolatory factorizations) into a factorization with
orthonormal basis functions. This process also diagonalizes all matrices
$B_{\sigma_{1}\sigma_{2}}$.

Suppose that the matrices $U_{\tau}$, $D_{\tau}$, and
$B_{\sigma_{1}\sigma_{2}}$ in an HSS factorization of a symmetric
matrices have already been generated (for instance by the algorithm
of Section \ref{sec:compHSS-ID}). The method described here produces
new matrices $U_{\tau}^{\rm new}$ and $B_{\sigma_{1}\sigma_{2}}^{\rm
new}$ with the property that each $U_{\tau}^{\rm new}$ has
orthonormal columns, and each $B_{\sigma_{1}\sigma_{2}}^{\rm new}$
is diagonal.

\begin{remark}
The method described here does not in any way rely on particular properties
of the interpolative decomposition. In fact, it works for any input matrices $U_{\tau}$,
$D_{\tau}$, $B_{\sigma_{1}\sigma_{2}}$ for which the factorization (\ref{eq:telescope}) holds.
\end{remark}

The orthonormalization procedure works hierarchically, starting at the finest
level and working upwards. At the finest level, it loops over all sibling
pairs $\{\sigma_{1},\,\sigma_{2}\}$. It orthonormalizes the basis matrices
$U_{\sigma_{1}}$ and $U_{\sigma_{2}}$ by computing their $QR$ factorizations,
$$
[W_{1},\,R_{1}] = \texttt{qr}(U_{\sigma_{1}})
\qquad\mbox{and}\qquad
[W_{2},\,R_{2}] = \texttt{qr}(U_{\sigma_{2}}),
$$
so that
$$
U_{\sigma_{1}} = W_{1}\,R_{1}
\qquad\mbox{and}\qquad
U_{\sigma_{2}} = W_{2}\,R_{2},
$$
and $W_{1}$ and $W_{2}$ have orthonormal columns.
The matrices $R_{1}$ and $R_{2}$ are then used to update the diagonal block
$B_{\sigma_{1}\sigma_{2}}$ to reflect the change in basis vectors,
$$
\tilde{B}_{12} = R_{1}\,B_{\sigma_{1}\sigma_{2}}\,R_{2}^{\rm t}.
$$
Then $\tilde{B}_{12}$ is diagonalized via a singular value decomposition,
$$
\tilde{B}_{12} = \tilde{W}_{1}\,B_{\sigma_{1}\sigma_{2}}^{\rm new}\,\tilde{W}_{2}.
$$
The new bases for $\sigma_{1}$ and $\sigma_{2}$ are constructed by updating
$W_{1}$ and $W_{2}$ to reflect the diagonalization of $\tilde{B}_{12}$,
$$
U_{\sigma_{1}}^{\rm new} = W_{1}\,\tilde{W}_{1},
\qquad\mbox{and}\qquad
U_{\sigma_{2}}^{\rm new} = W_{2}\,\tilde{W}_{2}.
$$
Finally, the basis vectors for the parent $\tau$ of $\sigma_{1}$ and $\sigma_{2}$
must be updated to reflect the change in bases at the finer level,
$$
\hat{U}_{\tau} \leftarrow \mtwo{\tilde{W}_{1}^{\rm t}\,R_{1}}{0}{0}{\tilde{W}_{2}^{\rm t}\,R_{2}}\,\hat{U}_{\tau}.
$$

Once the finest level has been processed, simply move up to the next
coarser one and proceed analogously. A complete description of the
recompression scheme is given with the caption ``Algorithm 2''.

\begin{figure*}
\fbox{
\begin{minipage}{\textwidth}
\begin{tabbing}
\hspace{15mm} \= \kill
\textit{Input:}  \> The matrices $\hat{U}_{\tau}$, $B_{\sigma_{1}\sigma_{2}}$, $D_{\tau}$ in an
                    HSS factorization of a symmetric matrix $A$.\\
\textit{Output:} \> Matrices $\hat{U}_{\tau}^{\rm new}$, $B_{\sigma_{1}\sigma_{2}}^{\rm new}$,
                    and $D_{\tau}$ that form an HSS  factorization of $A$ such that\\
                 \> all $\hat{U}_{\tau}^{\rm new}$ have orthonormal columns and all $B_{\sigma_{1}\sigma_{2}}^{\rm new}$
                    are diagonal.\\
                 \> (The matrices $D_{\tau}$ remain unchanged.)
\end{tabbing}

\rule{\textwidth}{1pt}

\begin{tabbing}
\hspace{8mm} \= \hspace{8mm} \= \hspace{8mm} \= \kill
Set $U_{\tau}^{\rm tmp} = U_{\tau}$ for all leaf nodes $\tau$.\\
\textbf{loop} over levels, finer to coarser, $p = P-1,\,P-2,\,\dots,\,0$\\
\> \textbf{loop} over all nodes $\tau$ on level $p$\\
\> \> Let $\sigma_{1}$ and $\sigma_{2}$ denote the two sons of $\tau$.\\
\> \> $[W_{1},\,R_{1}] = \texttt{qr}(\hat{U}_{\sigma_{1}}^{\rm tmp})$\\
\> \> $[W_{2},\,R_{2}] = \texttt{qr}(\hat{U}_{\sigma_{2}}^{\rm tmp})$\\
\> \> $[\tilde{W}_{1},\,B_{\sigma_{1}\sigma_{2}}^{\rm new},\,\tilde{W}_{2}] =
       \texttt{svd}(R_{1}\,B_{\sigma_{1}\sigma_{2}}\,R_{2}^{\rm t})$\\
\> \> $\hat{U}_{\sigma_{1}}^{\rm new} = W_{1}\,\tilde{W}_{1}$\\
\> \> $\hat{U}_{\sigma_{2}}^{\rm new} = W_{2}\,\tilde{W}_{2}$\\
\> \> $\hat{U}_{\tau}^{\rm tmp} = \mtwo{\tilde{W}_{1}^{\rm t}\,R_{1}}{0}{0}{\tilde{W}_{2}^{\rm t}\,R_{2}}
                        \hat{U}_{\tau}$\\
\> \textbf{end loop}\\
\textbf{end loop}
\end{tabbing}

\rule{\textwidth}{1pt}

\textit{Remark:} In practice, we let the matrices $U_{\tau}^{\rm tmp}$ and $U_{\tau}^{\rm new}$
simply overwrite $U_{\tau}$.
\end{minipage}}

\vspace{2mm}

\textsc{Algorithm 2:} Orthonormalizing an HSS factorization
\end{figure*}

\subsection{Non-symmetric matrices}
\label{sec:compHSSnonsym} The extension of Algorithms 1 and 2 to the
case of non-symmetric matrices is straight-forward. In Algorithm 1,
we construct a set of sample matrices $\{S_{\tau}\}$ with the
property that the columns of each $S_{\tau}$ span the column space
of the corresponding HSS row block $A^{\rm row}_{\tau}$. Since $A$
is in that case symmetric, the columns of $S_{\tau}$ automatically
span the row space of $A^{\rm col}_{\tau}$ as well. For
non-symmetric matrices, we need to construct different sample
matrices $S^{\rm row}_{\tau}$ and $S^{\rm col}_{\tau}$ whose columns
span the column space of $A^{\rm row}_{\tau}$ and the row space of
$A^{\rm col}_{\tau}$, respectively. Note that in practice, we work
only with the subsets of all these matrices formed by the respective
spanning rows and columns; in the non-symmetric case, these may be
different. The algorithm is described in full with the caption
``Algorithm 3''.

The generalization of the orthonormalization technique is entirely analogous.

\begin{figure*}
\fbox{
\begin{minipage}{\textwidth}
\begin{tabbing}
\hspace{15mm} \= \kill
\textit{Input:}  \> A fast means of computing matrix-vector products $x \mapsto A\,x$ and
                    $x\mapsto A^{\rm t}\,x$.\\
                 \> A method for computing individual entries of $A$.\\
                 \> The HSS-rank $k$ of $A$.\\
                 \> A partitioning $\{I_{(P,j)}\}_{j=1}^{2^{P}}$ of the interval $[1,\,2,\,\dots,\,N]$.\\
\textit{Output:} \> Matrices $\hat{U}_{\tau}$, $B_{\sigma_{1}\sigma_{2}}$, $D_{\tau}$ that form an
                    HSS factorization of $A$.
\end{tabbing}

\rule{\textwidth}{1pt}

\begin{tabbing}
\hspace{4mm} \= \hspace{4mm} \= \hspace{4mm} \= \hspace{68mm} \= \kill
Generate two $N\times (k+10)$ Gaussian random matrices $R^{\rm row}$ and $R^{\rm col}$.\\
Evaluate $S^{\rm row} = A^{\rm t}\,R^{\rm row}$ and $S^{\rm col} = A\,R^{\rm col}$ using the fast matrix-vector multiplier.\\
\textbf{loop} over levels, finer to coarser, $p = P:(-1):1$\\
\> \textbf{loop} over all nodes $\tau$ on level $p$\\
\> \> \textbf{if} $\tau$ is a leaf node then\\
\> \> \> $I_{\rm loc}^{\rm row} = I_{\tau}$ \>
         $I_{\rm loc}^{\rm col} = I_{\tau}$\\
\> \> \> $R_{\rm loc}^{\rm row} = R(I_{\tau},:)$ \>
         $R_{\rm loc}^{\rm col} = R(I_{\tau},:)$\\
\> \> \> $S_{\rm loc}^{\rm row} = S^{\rm row}(I_{\tau},:) - A(I_{\tau},I_{\tau})        \,R_{\rm loc}^{\rm row}$\>
         $S_{\rm loc}^{\rm col} = S^{\rm col}(I_{\tau},:) - A(I_{\tau},I_{\tau})^{\rm t}\,R_{\rm loc}^{\rm col}$\\
\> \> \textbf{else}\\
\> \> \> Let $\sigma_{1}$ and $\sigma_{2}$ be the two sons of $\tau$.\\
\> \> \> $I_{\rm loc}^{\rm row} = [\tilde{I}_{\sigma_{1}}^{\rm row},\,\tilde{I}_{\sigma_{2}}^{\rm row}]$\>
         $I_{\rm loc}^{\rm col} = [\tilde{I}_{\sigma_{1}}^{\rm col},\,\tilde{I}_{\sigma_{2}}^{\rm col}]$\\
\> \> \> $R_{\rm loc}^{\rm row} = \vtwo{R_{\sigma_{1}}^{\rm row}}{R_{\sigma_{2}}^{\rm row}}$\>
         $R_{\rm loc}^{\rm col} = \vtwo{R_{\sigma_{1}}^{\rm col}}{R_{\sigma_{2}}^{\rm col}}$\\
\> \> \> $S_{\rm loc}^{\rm row} = \vtwo{S_{\sigma_{1}}^{\rm row} - A(\tilde{I}_{\sigma_{1}}^{\rm row},\tilde{I}_{\sigma_{2}}^{\rm col})\,R_{\sigma_{2}}^{\rm row}}
                                       {S_{\sigma_{2}}^{\rm row} - A(\tilde{I}_{\sigma_{2}}^{\rm row},\tilde{I}_{\sigma_{1}}^{\rm col})\,R_{\sigma_{1}}^{\rm row}}$\>
         $S_{\rm loc}^{\rm col} = \vtwo{S_{\sigma_{1}}^{\rm col} - A(\tilde{I}_{\sigma_{1}}^{\rm row},\tilde{I}_{\sigma_{2}}^{\rm col})\,R_{\sigma_{2}}^{\rm col}}
                                       {S_{\sigma_{2}}^{\rm col} - A(\tilde{I}_{\sigma_{2}}^{\rm row},\tilde{I}_{\sigma_{1}}^{\rm col})\,R_{\sigma_{1}}^{\rm col}}$\\
\> \> \textbf{end if}\\
\> \> $[\hat{X}_{\tau}^{\rm row},J_{\tau}^{\rm row}] = \texttt{interpolate}((S_{\rm loc}^{\rm row})^{\rm t})$\>\>
      $[\hat{X}_{\tau}^{\rm col},J_{\tau}^{\rm col}] = \texttt{interpolate}((S_{\rm loc}^{\rm col})^{\rm t})$\\
\> \> $R_{\tau}^{\rm row} = (\hat{X}_{\tau}^{\rm col})^{\rm t}\,R_{\rm loc}^{\rm row}$\> \>
      $R_{\tau}^{\rm col} = (\hat{X}_{\tau}^{\rm row})^{\rm t}\,R_{\rm loc}^{\rm col}$\\
\> \> $S_{\tau}^{\rm row} = S_{\rm loc}^{\rm row}(J_{\tau}^{\rm row},:)$\> \>
      $S_{\tau}^{\rm col} = S_{\rm loc}^{\rm col}(J_{\tau}^{\rm col},:)$\\
\> \> $\tilde{I}_{\tau}^{\rm row} = I_{\rm loc}^{\rm row}(J_{\tau}^{\rm row})$\> \>
      $\tilde{I}_{\tau}^{\rm col} = I_{\rm loc}^{\rm col}(J_{\tau}^{\rm col})$\\
\> \textbf{end loop}\\
\textbf{end loop}\\
For all leaf nodes $\tau$, set $D_{\tau} = A(I_{\tau},\,I_{\tau})$.\\
For all sibling pairs $\{\sigma_{1},\,\sigma_{2}\}$ set
$B_{\sigma_{1}\sigma_{2}} = A(\tilde{I}_{\sigma_{1}}^{\rm row},\,\tilde{I}_{\sigma_{2}}^{\rm col})$.
\end{tabbing}
\end{minipage}}

\vspace{2mm}

\textsc{Algorithm 3:} Computing the HSS factorization of a non-symmetric matrix.

\end{figure*}

%
%


\providecommand{\bysame}{\leavevmode\hbox to3em{\hrulefill}\thinspace}
\providecommand{\MR}{\relax\ifhmode\unskip\space\fi MR }
\providecommand{\MRhref}[2]{%
  \href{http://www.ams.org/mathscinet-getitem?mr=#1}{#2}
}
\providecommand{\href}[2]{#2}

\end{document}